\documentclass[12pt,leqno]{amsart}
\usepackage{amssymb,amsthm}
\usepackage{amsmath,amssymb,color}
\numberwithin{equation}{section}

\newcommand{\ep}{\varepsilon}
\newcommand{\la}{\lambda}
\newcommand{\va}{\varphi}
\newcommand{\ppp}{\partial}

\newcommand{\www}{\widetilde}

\newcommand{\ddda}{d_t^{\alpha}}

\newcommand{\R}{\mathbb{R}}
\newcommand{\Q}{\mathbb{Q}}

\newcommand{\C}{\mathbb{C}}
\newcommand{\N}{\mathbb{N}}

\newcommand{\ooo}{\overline}

\newcommand{\OOO}{\Omega}
\newcommand{\MLONE}{E_{\alpha,1}}

\newcommand{\MLAA}{E_{\alpha,\alpha}}

\newcommand{\sumij}{\sum_{i,j=1}^d}

\newcommand{\hhalf}{\frac{1}{2}}
\newcommand{\DDD}{\mathcal{D}}

\newcommand{\sumn}{\sum_{n=1}^{\infty}}

\allowdisplaybreaks

\begin{document}
\title
[]
{
Simultaneous determination of initial value and source term
for time-fractional wave-diffusion equations}

\pagestyle{myheadings}
\markboth{P. Loreti, D. Sforza and M.~Yamamoto}
\author{$^1$ Paola Loreti, 
$^1$ Daniela Sforza,
$^{2,3,4}$ M.~Yamamoto }

\thanks{
$^1$  
Dipartimento di Scienze di Base e Applicate per l'Ingegneria, Sapienza 
Universit\`a di Roma, Via Antonio Scarpa 16, 00161 Roma, Italy 
e-mail: {\it loreti@dmmm.uniroma1.it}\\
$^1$
Dipartimento di Scienze di Base e Applicate per l'Ingegneria, Sapienza 
Universit\`a di Roma, Via Antonio Scarpa 16, 00161 Roma, Italy 
e-mail: {\it daniela.sforza@sbai.uniroma1.it}
\\
$^2$ Graduate School of Mathematical Sciences, The University
of Tokyo, Komaba, Meguro, Tokyo 153-8914, Japan \\
$^3$ Honorary Member of Academy of Romanian Scientists,
Ilfov, nr. 3, Bucuresti, Romania \\
$^4$ Correspondence member of Accademia Peloritana dei Pericolanti,\\
Palazzo Universit\`a, Piazza S. Pugliatti 1 98122 Messina Italy \\
e-mail: {\tt myama@ms.u-tokyo.ac.jp}
}

\date{}
\maketitle

\baselineskip 18pt

\begin{abstract}
We consider initial boundary value problems for
time fractional diffusion-wave equations:
$$
\ddda u = -Au + \mu(t)f(x)
$$
in a bounded domain where $\mu(t)f(x)$ describes a source and
$\alpha \in (0,1) \cup (1,2)$, and $-A$ is a symmetric 
ellitpic operator with repect to the spatial variable $x$.
We assume that $\mu(t) = 0$ for $t > T$:some time and 
choose $T_2>T_1>T$.
We prove the uniqueness in simultaneously determining $f$ in $\OOO$,
$\mu$ in $(0,T)$, and initial values of $u$ 
by data $u\vert_{\omega\times (T_1,T_2)}$, 
provided that the order $\alpha$ does not belong to 
a countably infinite set in $(0,1) \cup (1,2)$ which is characterized 
by $\mu$.
The proof is based on the asymptotic behavior of the Mittag-Leffler
functions.
\end{abstract}

\section{Introduction and Main Results}

Let $\OOO\subset \R^d$ be a bounded domain
with smooth boundary $\ppp\OOO$ and let $\nu = (\nu_1(x),..., 
\nu_d(x))$ be the unit outward normal vector to $\ppp\OOO$ at 
$x \in \ppp\OOO$. For $\alpha \in (0,1) \cup (1,2)$,
let $\ddda$ denote the Caputo derivative:
$$
\ddda v(t) = 
\left\{ \begin{array}{rl}
\frac{1}{\Gamma(1-\alpha)}\int^t_0 (t-s)^{-\alpha}
\frac{dv}{ds}(s) ds, \quad & 0<\alpha<1, \\
\frac{1}{\Gamma(2-\alpha)}\int^t_0 (t-s)^{1-\alpha}
\frac{d^2v}{ds^2}(s) ds, \quad & 1<\alpha<2
\end{array}\right.
$$
for suitable $v$, where $\Gamma(\cdot)$ denotes the gamma 
function.

We consider 
$$
\left\{ \begin{array}{rl}
& \ddda u = -Au + \mu(t)f(x), \quad x\in \OOO, t>0, \\
& u\vert_{\ppp\OOO\times (0,\infty)} = 0, \\
& u(\cdot, 0) = a \quad \mbox{in $\OOO$ $\quad$ if $0<\alpha < 1$}
\end{array}\right.
                          \eqno{(1.1)}
$$
and
$$
\left\{ \begin{array}{rl}
& \ddda u = -Au + \mu(t)f(x), \quad x\in \OOO, t>0, \\
& u\vert_{\ppp\OOO\times (0,\infty)} = 0, \\
& u(\cdot, 0) = a, \quad \ppp_tu(\cdot,0) = b \quad \mbox{in $\OOO$ 
$\quad$ if $1<\alpha < 2$}.
\end{array}\right.
                           \eqno{(1.2)}
$$
Here 
$$
-Au(x) = \sumij \ppp_i(a_{ij}(x)\ppp_ju) + c(x)u, \quad 
c\le 0 \quad \mbox{on $\ooo\OOO$},     \eqno{(1.3)}
$$
and $a, b \in L^2(\OOO)$.
We assume
$$
\mu \in L^{\infty}(0,T), \quad
\mu \not\equiv 0 \quad \mbox{in $(0,\infty)$}, \quad
\mu(t) = 0\quad \mbox{for $t\ge T$}.
                               \eqno{(1.4)}
$$
Then it is known that there exists a unique solution to each of
(1.1) and (1.2) in suitable classes.  More precisely,
we show 
\\
{\bf Lemma 1.}
\\
{\it 
In (1.1) and (1.2) 
we assume $\mu \in L^{\infty}(0,\infty)$ and $f \in L^2(\OOO)$.
\\
{\bf (i): $0<\alpha<1$:} 
In (1.1), we assume that $a \in L^2(\OOO)$.
To (1.1), there exists a unique solution $u \in C([0,T];L^2(\OOO))
\cap C((0,T]; H^2(\OOO) \cap H^1_0(\OOO))$ to (1.1).
\\
{\bf (ii): $1<\alpha<2$:} 
In (1.2) we assume $a, b \in L^2(\OOO)$.  
To (1.2), there exists a unique solution $u \in C([0,T];L^2(\OOO))
\cap C^1([0,T];\DDD(A^{-\frac{1}{\alpha}})) \cap
C((0,T];H^2(\OOO)\cap H^1_0(\OOO))$.
}
\\

Here $A^{-\frac{1}{\alpha}}$ denotes the fractional power of the operator 
$A$ (e.g., Pazy \cite{Pa}) defined by (1.3) with the domain 
$\DDD(A) = H^2(\OOO) \cap H^1_0(\OOO)$.
The space $\DDD(A^{-\frac{1}{\alpha}})$ is wider than $L^2(\OOO)$ and 
the topology is weaker than the one in $L^2(\OOO)$. 
The proof of the lemma is done by the eigenfunction expansions (2.1) stated 
below and as for the details, see Sakamoto and Yamamoto \cite{SY} for 
example.  We remark that the framework for coverning more general 
forms of equations with non-symmetric $A$ is available and see 
Kubica, Ryszewska and Yamamoto \cite{KRY}, Yamamoto \cite{Y2022} for
instance.
\\

Let $\omega \subset \OOO$ be an arbitrarily chosen subdomain 
and let $T<T_1<T_2$ be arbitrary.
We consider
\\
{\bf Inverse problem:} 
{\it Determine $f(x)$ and $\mu(t)$ of a source term and
initial values $a$ for $0<\alpha<1$ and $a, b$ for $1<\alpha<2$
simultaneously by data $u\vert_{\omega\times (T_1,T_2)}$.}
\\

For the statement of the main result, we introduce an index $\ell_0
\in \N \cup \{0\}$.  Here and henceforth let $\N:= \{1,2,3, ... \}$.

We show
\\
{\bf Lemma 2.}
\\
{\it
The space Span $\{ t^{\ell};\, \ell \in \N \cup \{0\}\}$ spanned by all the 
monomials is dense in 
$L^2(0,T)$.
}

The proof is direct from the Weierstrass polynomial approximation 
theorem and for completeness we prove the lemma at the end of Secton 2.

By $\mu\not\equiv 0$ in $(0,T)$, Lemma 2 yields 
that 
$$
\left\{ m\in \N \cup \{0\};\, \int^T_0 (-s)^m \mu(s) ds \ne 0
\right\} \ne \emptyset.
$$
Indeed, if not, then $\mu \in L^{\infty}(0,T) \subset L^2(0,T)$ is 
orthogonal to the dense set in $L^2(0,T)$, and so $\mu = 0$ in 
$(0,T)$, which is impossible by $\mu \not\equiv 0$ in $(0,T)$. 

Therefore, we can define 
$$
\ell_0:=  \min \left\{ m\in \N \cup \{0\};\, \int^T_0 (-s)^m \mu(s) ds 
\ne 0\right\}.
$$
In particular, we see $\ell_0=0$ if $\int^T_0 \mu(s) ds \ne 0$, and
if $\ell_0 \in \N$, then 
$$
\int^T_0 (-s)^{\ell_0} \mu(s) ds \ne 0, \quad
\int^T_0 (-s)^m \mu(s) ds = 0 \quad \mbox{for $0\le m \le \ell_0-1$}.
                                   \eqno{(1.5)}
$$
If $\mu \ge 0$ in $(0,T)$ satisfies (1.4), then $\int^T_0 \mu(s) ds 
\ne 0$, and so $\ell_0 = 0$.

Now we are ready to state the main result.
\\
{\bf Theorem 1.}
\\
{\it
{\bf (i): $0<\alpha<1$}.
We assume 
$$
\alpha \not\in \left\{ \frac{\ell_0+1}{n}\right\}_{n\in \N}.
                                          \eqno{(1.6)}
$$
Then $u\vert_{\omega\times (T_1,T_2)} = 0$ implies
$a=f=0$ in $\OOO$.
\\
{\bf (ii): $1<\alpha<2$}.
We assume 
$$
\alpha \not\in \left\{ \frac{\ell_0+1}{n}\right\}_{n\in \N}
\cup \left\{ \frac{\ell_0+ 2}{n}\right\}_{n\in \N}.
                                                    \eqno{(1.7)}
$$
Then $u\vert_{\omega\times (T_1,T_2)} = 0$ implies
$a=b=f=0$ in $\OOO$.
}
\\

In (1.6), (1,7) and also (1.11), (1.12) stated below, the 
right-hand sides mean sequences, that is, $\ell_0+1$ and $n$ are not assumed 
to be irreducible for example in (1.6).

For example, we see
\begin{align*}
& \left\{ \frac{\ell_0+1}{n} \right\}_{n\in \N}
= \left\{ \frac{1}{n} \right\}_{n\in \N} \quad \mbox{if
$\ell_0=0$}, \\
& \left\{ \frac{\ell_0+1}{n} \right\}_{n\in \N}
= \left\{ \frac{1}{n} \right\}_{n\in \N}
\cup \left\{ \frac{2}{2n-1} \right\}_{n\in \N}
\quad \mbox{if $\ell_0=1$}. \\
\end{align*}

We do not know whether we can remove the conditions (1.6) and (1.7)
on $\alpha$ in Theorem 1, and also (1.10) and (1.11) in Theorem 2 stated below.
We give more remark in 1. of Section 5. 

Now we can extend the uniqueness in Theorem 1 to the inverse problem of
determining both $f(x)$ and $\mu(t)$ as well as initial values.

For the statement of Theorem 2, in addition to (1.1) and (1.2), we consider 
$$
\left\{\begin{array}{rl}
& \ddda \www{u} = -A\www{u} + \www{\mu}(t)\www{f}(x), \quad x\in \OOO, t>0, \\
& \www{u}\vert_{\ppp\OOO\times (0,\infty)} = 0, \\
& \www{u}(\cdot, 0) = \www{a}, \quad \mbox{in $\OOO$ $\quad$  
if $0<\alpha < 1$}
\end{array}\right.
                           \eqno{(1.8)}
$$
and
$$
\left\{\begin{array}{rl}
& \ddda \www{u} = -A\www{u} + \www{\mu}(t)\www{f}(x), \quad x\in \OOO, t>0, \\
& \www{u}\vert_{\ppp\OOO\times (0,\infty)} = 0, \\
& \www{u}(\cdot, 0) = \www{a}, \quad \ppp_t\www{u}(\cdot,0) 
= \www{b} \quad \mbox{in $\OOO$ $\quad$
if $1<\alpha < 2$}.
\end{array}\right.
                           \eqno{(1.9)}
$$
We assume that $a, \,\www{a}, \,b, \,\www{b},\, f, \, \www{f}
\in L^2(\OOO)$ and $\mu, \www{\mu}$ satisfy (1.4).  Moreover 
let $f \not\equiv 0$ in $\OOO$ or $\www{f}\not\equiv 0$ in 
$\OOO$.
We define $\ell_1 \in \N \cup \{0\}$ as the minimal number satisfying 
$$
\ell_1 := \min\left\{ m\in \N \cup \{0\};\, 
\int^T_0 \mu(s)(-s)^m ds \ne 0 \quad \mbox{or} \quad 
\int^T_0 \www{\mu}(s)(-s)^m ds \ne 0\right\}.
$$
More precisely, we define $\ell_1 = 0$ if 
$\int^T_0 \mu(s) ds \ne 0$ or $\int^T_0 \www{\mu}(s)(-s)^m ds \ne 0$, and
$\ell_1 \ge 0$ if and only if 
$$
\left\{ \begin{array}{rl}
& \int^T_0 \mu(s)(-s)^m ds =  \int^T_0 \www{\mu}(s)(-s)^m ds = 0
\quad \mbox{for $0\le m \le \ell_1-1$},\\
& \int^T_0 \mu(s)(-s)^{\ell_1} ds \ne 0 \quad \mbox{or} \quad 
\int^T_0 \www{\mu}(s)(-s)^{\ell_1} ds \ne 0.
\end{array}\right.
                               \eqno{(1.10)}
$$
By $\mu\not\equiv 0$ and $\www{\mu} \not\equiv 0$, Lemma 2
guarantees the existence of such $\ell_1 \not\in \N \cup \{0\}$.
\\
{\bf Theorem 2.}
\\
{\it
{\bf (i): $0 < \alpha < 1$}.  Let $u, \www{u}$ satisfy (1.1) and (1.8) 
respectively.
We assume 
$$
\alpha \not\in \left\{ \frac{\ell_1+1}{n}\right\}_{n\in \N}.
                                                       \eqno{(1.11)}
$$
If $u = \www{u}$ in $\omega \times (T_1,T_2)$, then 
$a = \www{a}$ in $\OOO$ and there exists constants $c_1, c_2$ 
satisfying $\vert c_1\vert + \vert c_2 \vert > 0$ such that 
$$
c_1f(x) = c_2\www{f}(x), \quad x\in \OOO \quad \mbox{and}\quad
c_2\mu(t) = c_1\www{\mu}(t), \quad 0<t<T.
$$
\\
{\bf (ii): $1<\alpha<2$}.  Let $u, \www{u}$ satisfy (1.2) and (1.9) 
respectively.
We assume 
$$
\alpha \not\in \left\{ \frac{\ell_1+1}{n}\right\}_{n\in \N}
\cup \left\{ \frac{\ell_1+2}{n}\right\}_{n\in \N}
                                                       \eqno{(1.12)}
$$
If $u = \www{u}$ in $\omega \times (T_1,T_2)$, then 
$a = \www{a}$ and $b=\www{b}$ in $\OOO$ and there exists constants $c_1, c_2$ 
satisfying $\vert c_1\vert + \vert c_2 \vert > 0$ such that 
$$
c_1f(x) = c_2\www{f}(x), \quad x\in \OOO \quad \mbox{and}\quad
c_2\mu(t) = c_1\www{\mu}(t), \quad 0<t<T.
$$
}
\\

In (1.11) and (1.12), we recall that $\frac{\ell_1+1}{n}$ and 
$\frac{\ell_1+2}{n}$ are not assumed to be irreducible.

We remark that $\mu(t)f(x) = c\mu(t)\frac{1}{c}f(x)$ for $x\in \OOO$ and 
$t>0$ with any constant $c\ne 0$, we can at best expect the uniqueness 
up to a multiplicative constant, as is stated in Theorem 2.

Theorem 1 follows directly from Theorem 2 by setting $\www{\mu} \equiv 
\mu$ in $(0,T)$, and the proof is based on the same idea for 
Theorem 1.  For more readability, we separately state and prove
both theorems.

As for inverse problems of determining source terms of 
fractional differential equations in the case 
where extra data are taken along a time interval from $t=0$,
there are already many works.  Here not aiming at any comprehensive
list, we can refer to 
Jiang, Li, Liu and Yamamoto \cite{JLLY}, Jin, Kian and Zhou \cite{JKZ},
Kian \cite{K}, Kian, Soccorsi, Xue and Yamamoto \cite{KSXY},
Kian and Yamamoto \cite{KY}, Li, Liu and Yamamoto \cite{LLY1},
Liu, Li and Yamamoto \cite{LLY2}, Liu, Rundell and Yamamoto \cite{LRY},
Liu and Zhang \cite{LZ} and the references of these articles.
For determination of initial values, we here only refer to 
Jiang, Li, Pauron and Yamamoto \cite{JLPY},
Theorem 4.2 in Sakamoto and Yamamoto \cite{SY}.

We emphasize that unlike the above articles,
we exclusively take data not from 
$t=0$ but over time interval
$0< T_1 < t < T_2$.  It is more realistic than data starting at $t=0$ 
that we can start observations after $T>0$, that is, 
after the incident finishes at $t=T$, which is described by
$\mu(t) = 0$ for $t \ge T$.
For such a formulation of inverse source problems of determining $f(x)$ 
with known initial values, 
see Cheng, Lu and Yamamoto \cite{CLY} for 
parabolic equations and hyperbolic equations.
As for inverse problems of determining a source  
$f(x)$ with such data for time-fractional differential 
equations, see Jannno and Kian \cite{JaK}, Kian, Liu and Yamamoto
\cite{KLY}, Kinash and Janno \cite{KJ}, Yamamoto \cite{Ya}.

All the above referred works discuss two types for inverse source problems:
\begin{itemize}
\item
determination of source terms such as $f(x)$ with given initial values
$a$ and/or $b$.
\item
determination of initial values with given source terms.
\end{itemize}

To the best knowledge of the authors, there are no works on the 
simultaneous determination of $f(x)$ and $\mu(t)$ of a source
and initial values by single observation data $u\vert_{\omega\times (T_1,T_2)}$
after the incident finishes.
Thus Theorems 2 is the first uniqueness result for determination 
$f(x)$ and $\mu(x)$ as well as initial values.

The time-fractional diffusion-wave equations have remarkably 
different properties from the cases $\alpha=1$ and $\alpha=2$, and 
one distinguished character is less smoothing or averaging effect 
for solutions,
which implies that original profiles of initial values can be 
preserved more compared with the case $\alpha=1$, but 
less compared with $\alpha=2$.  Such a difference brings 
drastic differences especially for inverse problems for 
time-fractional diffusion-wave equations, and such an 
example is well-posedness of the backward problems in time, which 
is a strong contrast to $\alpha=1$, and we can refer only to 
Floridia, Li, Yamamoto \cite{FLY}, Floridia and Yamamoto \cite{FY}
as a few recent works.

Moreover by such less smoothing effect, we can expect the uniqueness in 
simultaneously detemining initial values and source terms, because 
the effect from initial values and the one from the source
can be separately reflected in the observation data.
This can not be correct for the cases $\alpha=1$ and even $\alpha=2$, as the 
following examples show for cases of ordinary differential equations.
\\
{\bf Example.} 
\\
Let $\la \ne 0$ be an arbitrarily chosen real constant. 
Let $a, f \in \R$ and 
$$
\mu(t) = \chi_{(0,T)}(t)
= \left\{ \begin{array}{rl}
& 1, \quad 0<t<T,\\
& 0, \quad t \ge T.
\end{array}\right.
$$
\\
{\bf (i) Case $\alpha=1$:}
We consider
$$
\frac{du}{dt}(t) = -\la u(t) + \chi_{(0,T)}(t)f, \quad u(0) = a.
$$ 
Then 
$$
u(t) = 
\left\{ \begin{array}{rl}
& e^{-\la t} a + e^{-\la t}\frac{e^{\la T} - 1}{\la}f, \quad 
t>T, \\
& e^{-\la t} a + \frac{1 - e^{-\la t}}{\la}f, \quad 0<t<T.
\end{array}\right.
$$
Therefore $u(t) = 0$ for $T_1<t<T_2$ if and only if
$a\la + (e^{\la T} - 1) f = 0$, which does not imply the 
uniqueness for two constants $a, f$.
\\
{\bf (ii) Case $\alpha=2$:}
We consider
$$
\frac{d^2u}{dt^2}(t) = -\la u(t) + \chi_{(0,T)}(t)f, \quad u(0) = a,
\quad \frac{du}{dt}(0) = b
$$ 
with $\la = r^2 > 0$.
Then 
\begin{align*}
& u(t) = a \cos rt + b\frac{\sin rt}{r}
+ f\int^T_0 \frac{\sin r(t-s)}{r} ds
= a \cos rt + b\frac{\sin rt}{r}
+ \frac{\cos r(T-t) - \cos rt}{r^2}f \\
=& \left( a + \frac{\cos rT -1}{r^2}f\right)\cos rt 
+ \left(\frac{b}{r} + \frac{\sin rT}{r^2} f\right) \sin rt, \quad
t>T.
\end{align*}
Therefore, $u(t) = 0$ for $T_1<t<T_2$ if and only if
$$
a + \frac{\cos rT -1}{r^2}f = \frac{b}{r} + \frac{\sin rT}{r^2} f = 0.
$$
In general we cannot conclude that $a=b=f=0$.
\\

This article is composed of five sections.
Section 2 provides some lemmata which are used for the proofs of 
Theorems 1 and 2, and Sections 3 and 4 are devoted to the proofs of Theorems 
1 and 2 respectively.
In Section 5, we give concluding remarks.

\section{Preliminaries}

It is known that the elliptic operator $A$ with the domain 
$\DDD(A) = H^2(\OOO) \cap H^1_0(\OOO)$ has a countably infinite number of 
eigenvalues, and by $c\le 0$ in $\OOO$ in (1.3), we can number as 
$$
0 < \la_1 < \la_2 < \cdots \longrightarrow \infty.
$$
Let $P_n$ be the eigenprojection for the eigenvalue $\la_n$ of 
$A$: $AP_nv = \la_n P_nv$ for each $n\in \N$.
We write  
$$
\Vert f\Vert:= \Vert f\Vert_{L^2(\OOO)}, \quad 
(f,g):= \int_{\OOO} f(x) g(x) dx.
$$

Henceforth by $E_{\beta, \gamma}(z)$ with $\beta>0$ and $\gamma > 0$,
we denote the Mittag-Leffler function (e.g., Podlubny \cite{Po}):
$$
E_{\beta, \gamma}(z):= \sum_{k=0}^{\infty}\frac{z^k}{\Gamma(\beta k + \gamma)},
\quad z\in \C.
$$
The power series can readily verify that $E_{\beta,\gamma}(z)$ is an 
entire function in $z\in \C$.

First we can show
\\
{\bf Lemma 3 (eigenfunction expansions)}.
\\
{\it
Let $u$ be the solutions to (1.1) and (1.2) in Lemma 1.
Then
$$
u(x,t) = 
\left\{ \begin{array}{rl}
& \sumn \MLONE(-\la_nt^{\alpha})P_na \\
+ &\sumn \left( \int^T_0 (t-s)^{\alpha-1}\MLAA(-\la_n(t-s)^{\alpha})
\mu(s) ds\right) P_nf, \quad t>T \quad \mbox{if $0<\alpha<1$}, \\
& \sumn \MLONE(-\la_nt^{\alpha})P_na 
+ \sumn tE_{\alpha,2}(-\la_nt^{\alpha}) P_nb \\
+ &\sumn \left(\int^T_0 (t-s)^{\alpha-1}\MLAA(-\la_n(t-s)^{\alpha})
\mu(s) ds\right) P_nf, \quad t>T \quad \mbox{if $1<\alpha<2$}.
\end{array}\right.
                                         \eqno{(2.1)}
$$
The series are convergent in the corresponding spaces indicated in 
Lemma 1.
}

The proof can be found for example, in \cite{SY}.

Next 
\\
{\bf Lemma 4.}
\\
{\it
The mapping $t \longrightarrow u(\cdot,t)$ is analytic from $(T,\infty)$ to 
$L^2(\OOO)$.
}
\\
{\bf Proof.}  The proof is done by the eigenfunction expansion (2.1).
\\
\vspace{0.2cm}

We further show
\\
{\bf Lemma 5.}
\\
{\it 
Let $N\in \N$ be fixed.
We assume that $p_1, ..., p_N$ are mutually distinct positive constants 
and $a_1, ..., a_N \in \R$.
If
$$
\sum_{k=1}^N \frac{a_k}{t^{p_k}} + O\left( \frac{1}{t^r}\right) = 0
\quad \mbox{as $t \to \infty$},
$$
then $a_k=0$ as long as $p_k<r$ with $k=1,2,..., N$.
}
\\
{\bf Proof of Lemma 5.}
Without loss of generality, we can assume 
that $p_1<p_2 < \cdots < p_N$.
We have
$$
a_1 + \sum_{k=2}^N a_k\frac{1}{t^{p_k-p_1}} = O\left( \frac{1}{t^{r-p_1}}
\right)
$$ 
as $t \to \infty$.  Therefore, $a_1=0$.
Next 
$$
\frac{a_2}{t^{p_2}} + \sum_{k=3}^N a_k\frac{1}{t^{p_k}} 
= O\left( \frac{1}{t^r} \right),
$$ 
that is,
$$
a_2 + \sum_{k=3}^N a_k\frac{1}{t^{p_k-p_2}} = O\left( \frac{1}{t^{r-p_2}}
\right)
$$ 
as $t \to \infty$.  Therefore, $a_2=0$.  We can continue this argument 
as long as $p_k < r$ for $k=1,2,..., N$.  Thus we 
completes the proof of Lemma 5.
$\blacksquare$

Furthermore, we show
\\
{\bf Lemma 6.}
\\
{\it 
Let 
$$
\sumn \vert a_n\vert^2 < \infty
$$
hold.  We assume that there exists $\{\ell_k\}_{k\in \N} 
\subset \N$ satisfying $\lim_{k\to \infty} \ell_k=\infty$.
If 
$$
\sumn \frac{a_n}{\la_n^{\ell_k}} = 0 \quad \mbox{for all $k\in \N$},
$$
then $a_n=0$ for all $n\in \N$.
}
\\
{\bf Proof of Lemma 6.}
We recall that $0<\la_1 < \la_2 < \cdots $ are the eigenvalues of $A$.
Setting $d_n:= \mbox{dim}\, \{w\in \DDD(A);\,
(A-\la_n)w = 0\}$, we create a sequence $\{ \mu_m\}_{m\in \N}$ by 
repeating $\la_n$ $d_n$-times, in other words, we number $\la_n$, $n\in \N$
according to the multiplicities $d_n$. 
Then it is known (e.g., Courant and Hilbert \cite{CH}) that 
$\mu_m \sim c_0m^{\frac{2}{d}}$ with some constant $c_0>0$ for all large 
$m\in N$.  By the definition of $\mu_m$, we have 
$\la_n \ge \mu_n$ for each $n\in \N$, and so 
$$
\la_n \ge c_0n^{\frac{2}{d}} \quad \mbox{for all large $n\in \N$}.
$$
Hence, since $\lim_{k\to\infty} \ell_k = \infty$, we see that there exists
$k_0\in \N$ such that 
$\frac{4\ell_k}{d} > 1$ for all $k\ge k_0$.
Consequently, by $\sumn \vert a_n\vert^2 < \infty$, we obtain
$$
\sumn \left\vert \frac{a_n}{\la_n^{\ell_k}} \right\vert
\le \left( \sumn \vert a_n\vert^2 \right)^{\hhalf}
\left( \sumn \left\vert \frac{1}{\la_n^{2\ell_k}}\right\vert\right)^{\hhalf}
\le C_1 \left(\sumn \vert a_n\vert^2\right)^{\hhalf}
\left(\sumn \frac{1}{n^{\frac{4\ell_k}{d}}}\right)^{\hhalf} < \infty
$$
for all $k\ge k_0$.

Since 
$$
\frac{a_1}{\la_1^{\ell_k}} + \sum_{n=2}^{\infty} 
\frac{a_n}{\la_n^{\ell_k}} = 0,\quad \mbox{that is},
\quad a_1 + \sum_{n=2}^{\infty} a_n \left( \frac{\la_1}{\la_n}\right)
^{\ell_k} = 0,
$$
we have
\begin{align*}
& \vert a_1\vert = \left\vert -\sum_{n=2}^{\infty} a_n 
\left( \frac{\la_1}{\la_n}\right)^{\ell_k}\right\vert
= \left\vert \sum_{n=2}^{\infty} 
\frac{a_n}{\la_n^{\ell_{k_0}}} \la_1^{\ell_{k_0}}
\left( \frac{\la_1}{\la_n}\right)^{\ell_k-\ell_{k_0}} \right\vert\\
\le & C\sum_{n=2}^{\infty} 
\left\vert \frac{a_n}{\la_n^{\ell_{k_0}}}\right\vert
\left( \frac{\la_1}{\la_n}\right)^{\ell_k-\ell_{k_0}} 
\le C_1\left( \frac{\la_1}{\la_2}\right)^{\ell_k-\ell_{k_0}}
\end{align*}
for all $k\ge k_0$.
By $0 < \la_1 < \la_2 < ....$, we see that 
$\left\vert \frac{\la_1}{\la_2}\right\vert < 1$.
Letting $k \to \infty$, we see that $\ell_k - \ell_{k_0} \to \infty$, and so
$a_1 = 0$.  Therefore, 
$$
\sum_{n=2}^{\infty} \frac{a_n}{\la_n^{\ell_k}} = 0, \quad k\in \N.
$$
Repeating the above argument, we have $a_2=a_3= \cdots = 0$.
Thus the proof of Lemma 6 is complete.
$\blacksquare$

We conclude this section with 
\\
{\bf Proof of Lemma 2.}
Since it is proved by the mollifier (e.g., Adams \cite{Ad}) 
that $C[0,T]$ is dense in $L^2(0,T)$, for any $u \in L^2(0,T)$ and 
$\ep > 0$, we can find $u_{\ep} \in C[0,T]$ such that 
$\Vert u - u_{\ep}\Vert_{L^2(0,T)} < \ep$.  By the 
Weierstrass polynomial approximation theorem, for any 
$\ep > 0$, we can find $p_{\ep} \in \, 
\mbox{Span}\, \{s^{\ell};\, \ell \in \N \cup \{ 0\}\}$ such that  
$\Vert p_{\ep} - u_{\ep}\Vert_{C[0,T]} < \ep$.  Therefore, 
\begin{align*}
& \Vert u - p_{\ep}\Vert_{L^2(0,T)} 
= \Vert u - u_{\ep} + u_{\ep} - p_{\ep} \Vert_{L^2(0,T)} 
\le \Vert u - u_{\ep}\Vert_{L^2(0,T)} 
 + \Vert u_{\ep} - p_{\ep} \Vert_{L^2(0,T)} \\
\le & \ep + \sqrt{T}\Vert u_{\ep} - p_{\ep} \Vert_{C[0,T]} 
\le (1+\sqrt{T})\ep.
\end{align*}
Thus the proof of Lemma 2 is complete.
$\blacksquare$

\section{Proof of Theorem 1}

The proof relies on the asymptotic expansions of the Mittag-Leffler functions.

We define a set $\{ m(k)\}_{k\in \N} \subset \N$ by 
$\{ k\in \N;\, \alpha k \not\in \N\}$.  
We see that $\{ m(k)\}_{k\in \N}$ is an infinite set.
\\
Indeed, if it is a finite set, then there exists some $N_0\in \N$ such that 
$n\alpha \in \N$ for $n \ge N_0$.
Therefore $N_0\alpha, \, (N_0+1)\alpha \in \N$, that is 
$\alpha = (N_0+1)\alpha - N_0\alpha \in \N$, which is impossible by 
$0<\alpha<1$ and $1<\alpha<2$.
$\blacksquare$
\\

More precisely, we can number
$$
m(k)\alpha \not\in \N, \quad m(1) < m(2) < \cdots \longrightarrow \infty.
$$
In particular, we note that $m(1) = 1$ and, if $\alpha \not\in \Q$, 
then $m(k) = k$ for all $k\in \N$, and that 
$$
\mbox{$j\in \N$ satisfies $\alpha j \not\in \N$ if and only if
$j=m(k)$ with some $k\in \N$.}
$$
\\
{\bf First Step: asymptotic expansions the solutions.}
\\
Let $\la > 0$.
The asymptotic behavior of the Mittag-Leffler functions (e.g., 
\cite{Po}) yields
$$
\MLONE(-\la t^{\alpha}) 
= \sum_{k\in \{1, ..., M\}, \alpha k\not\in \N}
\frac{(-1)^{k+1}}{\Gamma(1-\alpha k)}
\frac{1}{\la^k}\frac{1}{t^{\alpha k}}
+ O\left( \frac{1}{t^{\alpha(M+1)}}\right), 
$$
$$
t^{\alpha-1}\MLAA(-\la t^{\alpha}) 
= \sum_{k\in \{2,..., M\}, \alpha (k-1)\not\in \N}
 \frac{(-1)^{k+1}}{\Gamma((1-k)\alpha)}
\frac{1}{\la^k}\frac{1}{t^{\alpha (k-1)+1}}
+ O\left( \frac{1}{t^{\alpha M+1}}\right)
$$
and
$$
tE_{\alpha,2}(-\la t^{\alpha}) 
= \sum_{k\in \{1, ..., M\}, \alpha k\not\in \N\setminus \{1\}} 
\frac{(-1)^{k+1}}{\Gamma(2-\alpha k)}
\frac{1}{\la^k}\frac{1}{t^{\alpha k-1}}
+ O\left( \frac{1}{t^{\alpha(M+1)-1}}\right)
$$
for all large $t>0$ and all $M\in \N$.

Since $\frac{1}{\Gamma(-k)} = 0$ for $k \in \N \cup \{ 0\}$, 
we can represent the series as follows in terms of $m(k)$.
First
$$
E_{\alpha,1}(-\la t^{\alpha})
= \sum_{k=1}^N \frac{(-1)^{m(k)+1}}{\Gamma(1-\alpha m(k))}
\frac{1}{\la^{m(k)}} \frac{1}{t^{\alpha m(k)}}
+ O\left( t^{-\alpha m(N+1)}\right)
$$
as $t \to \infty$.

Next, for $j\in \N$, noting that $\alpha (j-1) \not\in \N$ if and only if
$j\in \{ m(k) + 1\}_{k\in \N}$, we represent 
$$
t^{\alpha-1}E_{\alpha,\alpha}(-\la t^{\alpha})
= \sum_{j=1}^N \frac{(-1)^{m(j)}}{\Gamma(-\alpha m(j))}
\frac{1}{\la^{m(j)+1}} \frac{1}{t^{\alpha m(j)+1}}
+ O\left( t^{-\alpha m(N+1)-1}\right)
$$
as $t \to \infty$.

Finally, for $1<\alpha<2$, we see that 
$\alpha k \not\in \N \setminus \{1\}$ if and only if 
$\alpha k \not\in \N$.
\\
Indeed, it is trivial that $\alpha k \not\in \N$ implies
$\alpha k \not\in \N \setminus \{1\}$.  Let
$\alpha k \not\in \N \setminus \{1\}$.  Then 
$\alpha k \not\in \N$ or $\alpha k = 1$ with $k\in \N$.
Then $\alpha k = 1$ implies $\alpha = \frac{1}{k} \le 1$, which contradicts 
that $\alpha > 1$.  Hence, $\alpha k \not\in \N\setminus \{ 1\}$ 
implies $\alpha k \not\in \N$.
$\blacksquare$
\\
Consequently, 
$$
t^{\alpha-1}E_{\alpha,\alpha}(-\la t^{\alpha})
= \sum_{j=1}^N \frac{(-1)^{m(j)+1}}{\Gamma(2-\alpha m(j))}
\frac{1}{\la^{m(j)}} \frac{1}{t^{\alpha m(j)-1}}
+ O\left( t^{-\alpha m(N+1)+1}\right)
$$
as $t \to \infty$ for each $N\in \N$.
\\

To sum up, 
$$
\left\{ \begin{array}{rl}
& \MLONE(-\la_n t^{\alpha}) 
= \sum_{k=1}^N  \frac{(-1)^{m(k)+1}}{\Gamma(1-\alpha m(k))}
\frac{1}{\la_n^{m(k)}}\frac{1}{t^{\alpha m(k)}}
+ O\left( t^{-\alpha m(N+1)}\right), \cr \\
& tE_{\alpha,2}(-\la_n t^{\alpha}) 
= \sum_{k=1}^N  \frac{(-1)^{m(k)+1}}{\Gamma(2-\alpha m(k))}
\frac{1}{\la_n^{m(k)}}\frac{1}{t^{\alpha m(k)-1}}
+ O\left( t^{-\alpha m(N+1)+1} \right) \quad 
                    \mbox{for $1<\alpha<2$},\cr \\
& t^{\alpha-1}\MLAA(-\la_n t^{\alpha}) 
= \sum_{k=1}^N \frac{(-1)^{m(k)}}{\Gamma(-\alpha m(k))}
\frac{1}{\la_n^{m(k)+1}}\frac{1}{t^{\alpha m(k)+1}}
+ O\left( t^{-\alpha m(N+1)-1}\right)
\end{array}\right.
                                     \eqno{(3.1)}
$$
for all large $t>0$ and each $n, N \in \N$.
\\ 
Consequently,
\begin{align*}
& \sumn (P_na)(x) \MLONE(-\la_n t^{\alpha}) \\
=& \sumn (P_na) \sum_{k=1}^N  \frac{(-1)^{m(k)+1}}{\Gamma(1-\alpha m(k))}
\frac{1}{\la_n^{m(k)}}\frac{1}{t^{\alpha m(k)}}
+ \sumn (P_na) O\left( t^{-\alpha m(N+1)}\right) \\
=& \sum_{k=1}^N  \frac{(-1)^{m(k)+1}}{\Gamma(1-\alpha m(k))}
\left( \sumn \frac{P_na}{\la_n^{m(k)}}\right) \frac{1}{t^{\alpha m(k)}}
+ O\left( t^{-\alpha m(N+1)}\right), \cr \\
& \sumn (P_nb)(x) tE_{\alpha,2}(-\la_n t^{\alpha}) \\
=& \sumn (P_nb) \sum_{k=1}^N  \frac{(-1)^{m(k)+1}}{\Gamma(2-\alpha m(k))}
\frac{1}{\la_n^{m(k)}}\frac{1}{t^{\alpha m(k)-1}}
+ \sumn (P_nb) O\left( t^{-\alpha m(N+1)+1}\right)   \\
=& \sum_{k=1}^N  \frac{(-1)^{m(k)+1}}{\Gamma(2-\alpha m(k))}
\left( \sumn \frac{P_nb}{\la_n^{m(k)}}\right) \frac{1}{t^{\alpha m(k)-1}}
+ O\left( t^{-\alpha m(N+1)+1}\right),   \cr \\
&\sumn \left( \int^T_0 (t-s)^{\alpha-1}\MLAA(-\la_n(t-s)^{\alpha})
\mu(s) ds \right)P_nf \\
=& \sumn \left(  \int^T_0 \sum_{k=1}^N  \frac{(-1)^{m(k)}}
{\Gamma(-\alpha m(k))}(t-s)^{-\alpha m(k) -1} \mu(s) ds\right)
\frac{1}{\la_n^{m(k)+1}}(P_nf)                                          \\
+ & \sumn (P_nf)\left( \int^T_0 
O((t-s)^{-\alpha m(N+1)-1}) \mu(s) ds\right) 
\end{align*}
$$
= \sum_{k=1}^N  \frac{(-1)^{m(k)}}
{\Gamma(-\alpha m(k))} \left( \int^T_0 (t-s)^{-\alpha m(k) -1} \mu(s) ds\right)
\sumn \frac{P_nf}{\la_n^{m(k)+1}}
+ O( t^{-\alpha m(N+1)-1})             \eqno{(3.2)}
$$
for all large $t>0$.

For the last equality, we used $(t-s)^{-\alpha m(N+1)-1}
\sim t^{-\alpha m(N+1)-1}$ for $t > 2T$ and $0 < s < T$, and 
$\int^T_0 O((t-s)^{-\alpha m(N+1)-1}) \mu(s) ds
= O(t^{-\alpha m(N+1)-1})$ as $t\to \infty$.

Henceforth we assume $t>2T$ and we set 
$$
\mu_m: = \int^T_0 (-s)^m \mu(s) ds, \quad m \in \N \cup \{0\}.
$$   
\\
{\bf Calculations of $\int^T_0 (t-s)^{-\sigma} \mu(s) ds$ with 
$\sigma = \alpha m(k) + 1$.}
\\
Let $t>2T$ and $\sigma>0$.  Then the binomial expansion implies
$$
(t-s)^{-\sigma} = \sum_{\ell=0}^L 
\left(
\begin{array}{c}
-\sigma \\
\ell \\
\end{array}
\right) t^{-\sigma-\ell}(-s)^{\ell} + 
R_L(t,s), \quad 0<s<T, \, t>2T,                                  \eqno{(3.3)}
$$
where we write 
$$
\left(
\begin{array}{c}
-\sigma \\
0 \\
\end{array}
\right) := 1, \quad
\left(
\begin{array}{c}
-\sigma \\
\ell \\
\end{array}
\right) := \frac{(-\sigma)(-\sigma-1) \cdots (-\sigma-\ell+1)}{\ell!}
\quad \mbox{for $\ell \in \N$},
$$
and the function $R_L(t,s)$ satisfies 
$$
\vert R_L(t,s)\vert \le C(L,\sigma)t^{-\sigma-L-1}T^{L+1}.
$$
Indeed for fixed $0<\rho<1$, noting that 
$(t-s)^{-\sigma} = t^{-\sigma}\left( 1 - \left(\frac{s}{t}\right)\right)
^{-\sigma}$, we can derive (3.2) from the Taylor series:
$$
(1-\eta)^{-\sigma} = \sum_{\ell=0}^L 
\left(
\begin{array}{c}
-\sigma \\
\ell \\
\end{array}
\right) (-\eta)^{\ell} + R_L(\eta), \quad \vert \eta\vert \le \rho,
$$
where 
$$
\vert R_L(\eta)\vert \le C(L,\sigma)\vert \eta\vert^{L+1} \quad 
\mbox{for $\vert \eta\vert \le \rho$}.
$$
Thus (3.3) is seen.
$\blacksquare$

We apply (3.3) to obtain
$$
\left\{ \begin{array}{rl}
& \int^T_0 (t-s)^{-\sigma} \mu(s) ds
= \left(
\begin{array}{c}
-\sigma \\
\ell_0 \\
\end{array}
\right) t^{-\sigma-\ell_0}\mu_{\ell_0} + O(t^{-\sigma-\ell_0-1}), 
\quad t>2T, \\
& \mbox{where $\ell_0=0$ or $\ell_0 \in \N$ satisfying 
$\int^T_0 (-s)^m \mu(s) ds = 0$ for $0\le m\le \ell_0-1$}.
\end{array}\right.
                                     \eqno{(3.4)}
$$
We note that 
$\left(
\begin{array}{c}
-\sigma \\
\ell_0 \\
\end{array}
\right) 
= \left(
\begin{array}{c}
-\alpha m(k) - 1 \\
\ell_0 \\
\end{array}
\right) \ne 0$ for each $k\in \N$, because 
$\alpha m(k) \not\in \N$ by the definition
implies $-\alpha m(k) - 1 \not\in \N$.
\\
{\bf Verification of (3.4).}
By (3.3) with $L=\ell_0$, since $\mu_0 = \cdots = \mu_{\ell_0-1} = 0$, 
we have
\begin{align*}
& \int^T_0 \mu(s) (t-s)^{-\sigma} ds\\
=& \sum_{\ell=0}^{\ell_0} \left(
\begin{array}{c}
-\sigma \\
\ell \\
\end{array}
\right) t^{-\sigma-\ell}\int^T_0 (-s)^{\ell} \mu(s) ds
+ \int^T_0 R_{\ell_0}(t,s) \mu(s) ds\\
= & \left(
\begin{array}{c}
-\sigma \\
\ell_0 \\
\end{array}
\right) t^{-\sigma-\ell_0}\mu_{\ell_0}
+ \int^T_0 R_{\ell_0}(t,s) \mu(s) ds.
\end{align*}
Moreover, 
$$
\left\vert \int^T_0 R_{\ell_0}(t,s) \mu(s) ds\right\vert
\le C\Vert \mu\Vert_{L^{\infty}(0,T)}t^{-\sigma-\ell_0-1}
T^{\ell_0+1}.
$$
Thus the verification of (3.4) is complete.
$\blacksquare$.

We set 
$$
\left\{ \begin{array}{rl}
& Q_k(x):= \frac{(-1)^{m(k)+1}}{\Gamma(1-\alpha m(k))}
\sumn \frac{P_na(x)}{\la_n^{m(k)}}, \cr\\
& R_k(x) := \frac{(-1)^{m(k)+1}}{\Gamma(2-\alpha m(k))}
\sumn \frac{P_nb(x)}{\la_n^{m(k)}}, \cr \\
& S_k(x):=  \frac{(-1)^{m(k)}}{\Gamma(-\alpha m(k))}
\left( \sumn \frac{(P_nf)(x)}{\la_n^{m(k)+1}}\right) 
\left(
\begin{array}{c}
-\alpha m(k) -1 \\
\ell_0 \\
\end{array}
\right).
\end{array}\right.
$$
Now we assume 
$$
u(x,t) = 0, \quad x\in \omega, \, T_1<t<T_2.
$$
Then Lemma 4 implies that $u(x,t)=0$ for $x\in \omega$ and 
all $t > T$.  

We note that $\mu_{\ell_0} \ne 0$ and 
$$
t^{-\alpha m(N+1) - 1} \le Ct^{-\alpha m(N+1)}, \quad
t^{-\alpha m(N+1) - 1} \le Ct^{-\alpha m(N+1)+1}
$$
for all $t > 2T$, and in terms of 
$\mu_{\ell_0}\left(
\begin{array}{c}
-\alpha m(k) -1 \\
\ell_0 \\
\end{array}
\right) \ne 0$, we obtain
\begin{align*}
& \int^T_0 (t-s)^{-\alpha m(k) - 1} \mu(s) ds
= \left(
\begin{array}{c}
-\alpha m(k) - 1 \\
\ell_0 \\
\end{array}
\right) t^{-\alpha m(k) - \ell_0-1}\mu_{\ell_0} 
+ O(t^{-\alpha m(k) - \ell_0-2})\\
=& t^{-\alpha m(k) - \ell_0-1}\mu_{\ell_0} 
\left(
\begin{array}{c}
-\alpha m(k) -1 \\
\ell_0 \\
\end{array}
\right) \left( 1 + O\left(\frac{1}{t}\right)\right).
\end{align*}
Substituting the 
above asymptotic expansions into (3.2) and (2.1), we obtain:
\\
{\bf Case: $0<\alpha<1$.}
$$
u(x,t) = \sum_{k=1}^N Q_kt^{-\alpha m(k)}
+ \sum_{k=1}^N S_k\mu_{\ell_0}
t^{-\alpha m(k) - \ell_0-1} 
$$
$$
+ O(t^{-\alpha m(1) - \ell_0 - 2})
+ O(t^{-\alpha m(N+1)}) = 0, \quad x\in \omega, \, t> 2T
                                               \eqno{(3.5)}
$$
for $N\in \N$.
\\
{\bf Case: $1<\alpha<2$.}
$$
u(x,t) = \sum_{k=1}^N Q_kt^{-\alpha m(k)}
+ \sum_{k=1}^N R_kt^{-\alpha m(k)+1}
$$
$$
+ \sum_{k=1}^N S_k\mu_{\ell_0}t^{-\alpha m(k) - \ell_0-1}
+ O(t^{-\alpha m(1) - \ell_0 - 2}) + O(t^{-\alpha m(N+1)+1}) = 0,
\quad x\in \omega, \, t>2T
                                               \eqno{(3.6)}
$$
for $N\in \N$.
\\
{\bf Second Step: Completion of the proof for $0<\alpha<1$.}
\\
Since $\alpha \not\in \left\{ \frac{\ell_0+1}{n}\right\}_{n\in \N}$, we
see 
$$
\{ \alpha m(k)\}_{k\in \N} \cap \{ \alpha m(k) + \ell_0 + 1\}_{k\in \N}
= \emptyset.                         \eqno{(3.7)}
$$
Indeed, if $\alpha m(k) = \alpha m(k') + \ell_0 + 1$ with some
$k, k' \in \N$, we have
$\alpha(m(k) - m(k')) = \ell_0 + 1$, that is,
$\alpha = \frac{\ell_0+1}{m(k) - m(k')}$, which is impossible by
(1.6): 
$\alpha \not\in \left\{ \frac{\ell_0 + 1}{n} \right\}_{n\in \N}$.
This is the verification of (3.7).
$\blacksquare$
\\

By (3.7), we see that $\alpha m(j), \alpha(k) + \ell_0 + 1$ for 
$1\le j,k \le N$ are mutually distinct.  We choose $N_1 \in \N$ 
large, so that 
$$
\alpha m(1), \alpha m(1) + \ell_0 + 1 < \alpha m(N_1+1).  
$$
Then 
$$
\alpha m(1), \alpha m(1) + \ell_0 + 1 
< \min\{ \alpha m(N_1+1), \, \alpha m(1) + \ell_0 + 2\} =: r,
$$
and (3.5) implies
$$
\sum_{k=1}^{N_1} Q_k t^{-\alpha m(k)}
+ \sum_{k=1}^{N_1} \mu_{\ell_0}S_k t^{-\alpha m(k)-\ell_0-1}
= O(t^{-r}) \quad \mbox{in $\omega$}     \eqno{(3.8)}
$$
as $t \to \infty$.   Hence we can apply Lemma 5 to (3.8), so that 
$Q_1 = S_1 = 0$ in $\omega$.

Therefore, in view of (3.5), we obtain
$$
\sum_{k=2}^N Q_k t^{-\alpha m(k)}
+ \sum_{k=2}^N S_k \mu_{\ell_0}t^{-\alpha m(k) - \ell_0-1}
+ O(t^{-\alpha m(2) - \ell_0 - 2})
= O(t^{-\alpha m(N+1)}) \quad \mbox{in $\omega$}
                                           \eqno{(3.9)}
$$
as $t \to \infty$, where $N\in \N$ is arbitrarily given.
We can repeat the previous argument to (3.9), and so 
we can obtain $Q_2 = S_2 = 0$ in $\omega$.
Thus, continuing this argument, we can complete the proof of 
$Q_1 = Q_2 = \cdots = Q_N = S_1 = S_2 = \cdots = S_N = 0$ in 
$\omega$ for each $N\in \N$.

Hence,
$$
\sumn \frac{P_na}{\la_n^{m(k)}} = \sumn \frac{P_nf}{\la_n^{m(k)+1}}
= 0 \quad \mbox{in $\omega$ for all $k \in \N$}.
$$
For large $k\in \N$, similarly to the proof of Lemma 6,
by using $\la_n \ge c_0n^{\frac{2}{d}}$ for large $n\in \N$,
we see that the series are convergent in $L^2(\omega)$. 
Therefore, for any $\va \in L^2(\omega)$, we see 
$$
\sumn \frac{(P_na,\va)_{L^2(\omega)}}{\la_n^{m(k)}}
= \sumn \frac{(P_nf,\va)_{L^2(\omega)}}{\la_n^{m(k)+1}} = 0
$$
for all $k\in \N$. Hence, Lemma 6 yields
$$
(P_na,\va)_{L^2(\omega)} = (P_nf,\va)_{L^2(\omega)} = 0 \quad
\mbox{in $\omega$ for all $n\in \N$ and all $\va \in L^2(\omega)$.}
$$
This means that $P_na = P_nf = 0$ in $\omega$ for all $n\in \N$.
Then we can prove that $P_na = P_nf = 0$ in $\OOO$ for all $n\in \N$.
\\
Indeed, the definition of $P_n$ implies 
$(A-\la_n)P_na = 0$ in $\OOO$. By $P_na = 0$ in $\omega$, the 
unique continuation (e.g., Isakov \cite{Is}) for the elliptic 
operator $A-\la_n$ yields that $P_na = 0$ in $\OOO$.
Similarly we can obtain $P_nf= 0$ in $\OOO$ for all $n\in \N$.
Since $a = \sumn P_na$ and $f = \sumn P_nf$ in $L^2(\OOO)$, we reach 
$a=f=0$ in $\OOO$.  Thus the proof of Theorem 1 is complete for the case 
$0<\alpha<1$.
$\blacksquare$
\\
\vspace{0.2cm}
\\
{\bf Third Step: Completion of the proof for $1<\alpha<2$.}
\\
Since $\alpha \not\in \left\{ \frac{\ell_0+1}{n}\right\}_{n\in \N}
\cup \left\{ \frac{\ell_0+2}{n}\right\}_{n\in \N}$ and 
$\alpha > 1$, we can verify
$$
\left\{ \begin{array}{rl}
& \{ \alpha m(k)\}_{k\in \N} \cup \{ \alpha m(k) - 1\}_{k\in \N}
= \emptyset,\\
& \{ \alpha m(k)\}_{k\in \N} \cup \{ \alpha m(k) + \ell_0 + 1\}_{k\in \N}
= \emptyset,\\
& \{ \alpha m(k)-1\}_{k\in \N} \cup \{ \alpha m(k) +\ell_0+1\}_{k\in \N}
= \emptyset.
\end{array}\right.
                                             \eqno{(3.10)}
$$
By (3.10), we see that 
$$
\alpha m(i), \, \alpha m(j)-1, \, \alpha m(k) + \ell_0 + 1, \quad
i,j,k \in \{ 1, 2, ..., N \}
$$
are mutually distinct.  We choose $N_1 \in \N$ such that 
$$
\alpha m(1), \, \alpha m(1) + \ell_0 + 1 < \alpha m(N_1+1) - 1.
                                               \eqno{(3.11)}
$$
Then $\alpha m(1)-1 < \alpha m(N_1+1) - 1$.

With this $N_1$, we can write (3.6) as 
$$
\sum_{k=1}^{N_1} Q_k t^{-\alpha m(k)} + \sum_{k=1}^{N_1}R_kt^{-\alpha m(k) + 1}
+ \sum_{k=1}^{N_1} S_k t^{-\alpha m(k) - \ell_0 - 1}
$$
$$
+ O(t^{-\alpha m(1) - \ell_0 - 2}) = O(t^{-\alpha m(N_1+1) + 1})                                                            \eqno{(3.12)}
$$
in $\omega$ as $t \to \infty$.  

In view of (3.11), we apply Lemma 5 to (3.12), and so we can obtain
$$
Q_1 = R_1 = S_1 = 0 \quad \mbox{in $\omega$}.    \eqno{(3.13)}
$$

Hence by (3.13), for all $N \in \N$ we obtain
$$
\sum_{k=s}^N Q_k t^{-\alpha m(k)} + \sum_{k=2}^N R_kt^{-\alpha m(k) + 1}
+ \sum_{k=2}^N S_k \mu_{\ell_0}t^{-\alpha m(k) - \ell_0 - 1}
$$
$$
+ O(t^{-\alpha m(2) - \ell_0 - 2}) 
= O(t^{-\alpha m(N+1) + 1})                  \eqno{(3.14)}
$$
in $\omega$ as $t \to \infty$.  
We apply the same argument as (3.13) with $N=N_1$ to (3.14), and so 
we obtain
$Q_2 = R_2 = S_2= 0$ in $\omega$.

Continuing this argument, we reach 
$Q_k = R_k = S_k = 0$ in $\omega$ for each $k\in \N$.
Similarly to the case $0<\alpha<1$, we can verify that 
$Q_k=R_k=S_k=0$ in $\omega$ for all $k\in \N$ imply
$P_na = P_nb = P_nf = 0$ in $\OOO$ for all $n\in\N$.
Therefore $a = b = f = 0$ in $\OOO$.
Thus the proof of Theorem 1 is complete also for $1<\alpha<2$.
$\blacksquare$

\section{Proof of Theorem 2}

The proof is the repeat of the same arguments as in the proof of 
Theorem 1 which is by equating the coefficients of the asymptotic 
expansions of the solutions.

For $\ell \in \N \cup \{ 0\}$, $k\in \N$ and $x\in \OOO$ and 
$a, b, f \in L^2(\OOO)$, we 
set  
$$
\left\{ \begin{array}{rl}
& Q_k(a)(x):= \frac{(-1)^{m(k)+1}}{\Gamma(1-\alpha m(k))}
\sumn \frac{P_na(x)}{\la_n^{m(k)}},                 \cr \\
& R_k(b)(x) := \frac{(-1)^{m(k)+1}}{\Gamma(2-\alpha m(k))}
\sumn \frac{P_nb(x)}{\la_n^{m(k)}}, \cr \\
& S_{k,\ell}(f)(x):=  \frac{(-1)^{m(k)}}{\Gamma(-\alpha m(k))}
\left(
\begin{array}{c}
-\alpha m(k) -1 \\
\ell \\
\end{array}
\right)
\left( \sumn \frac{(P_nf)(x)}{\la_n^{m(k)+1}}\right) 
\end{array}\right.
$$
and
$$
\mu_{\ell}:= \int^T_0 (-s)^{\ell} \mu(s) ds, \quad
\www{\mu}_{\ell}:= \int^T_0 (-s)^{\ell} \www{\mu}(s) ds \quad
\mbox{for all $\ell \in \N \cup \{ 0\}$}.
$$
Henceforth, without loss of generality, in (1.10), we assume that 
$\mu_{\ell_1} = \int^T_0 (-s)^{\ell_1} \mu(s) ds \ne 0$.
\\
{\bf Proof in the case $0<\alpha<1$.}
\\
Similarly to (3.5), by means of $u(x,t) = \www{u}(x,t)$ for 
$x\in \omega$ and $T_1<t<T_2$ and Lemma 4, we have
\begin{align*}
& \sum_{k=1}^N Q_k(a)t^{-\alpha m(k)}
+ \sum_{k=1}^N \mu_{\ell_1} S_{k,\ell_1}(f) t^{-\alpha m(k) - \ell_1-1}
=& \sum_{k=1}^N Q_k(\www{a})t^{-\alpha m(k)}
+ \sum_{k=1}^N \www{\mu}_{\ell_1} S_{k,\ell_1}(\www{f}) 
t^{-\alpha m(k) - \ell_1-1}\\
+ & O(t^{-\alpha m(1) - \ell_1-2}) + O(t^{-\alpha m(N+1)})
\end{align*}
in $\omega$ as $t \to \infty$.
Here we applied also 
$$
\int^T_0 (t-s)^{-\alpha m(k) -1} \mu(s) ds
= \mu_{\ell_1} 
\left(
\begin{array}{c}
-\alpha m(k) -1 \\
\ell_1 \\
\end{array}
\right)
t^{-\alpha m(k) - 1 - \ell_1}
+ O(t^{-\alpha m(k) - 2 - \ell_1})
$$
and
$$
\int^T_0 (t-s)^{-\alpha m(k) -1} \www{\mu}(s) ds
= \www{\mu}_{\ell_1} 
\left(
\begin{array}{c}
-\alpha m(k) -1 \\
\ell_1 \\
\end{array}
\right)
t^{-\alpha m(k) - 1 - \ell_1}
+ O(t^{-\alpha m(k) - 2 - \ell_1}),
$$
which are proved in the same way as for (3.3) by using
$\int^T_0 \mu(s) (-s)^m ds =  \int^T_0 \www{\mu}(s) (-s)^m ds
= 0$ for $0 \le m \le \ell_1 - 1$ with $\ell_1 \in \N$.
We note 
$$
\mu_{\ell_1} \ne 0,                        \eqno{(4.1)}
$$
but $\www{\mu}_{\ell_1} = 0$ may occur.
Since the definition of $S_{k,\ell}(\cdot)$ yields
$$
\mu_{\ell_1}S_{k,\ell_1}(f) - \www{\mu}_{\ell_1}S_{k,\ell_1}(\www{f})
= S_{k,\ell_1}(\mu_{\ell_1}f - \www{\mu}_{\ell_1}\www{f})
\quad \mbox{in $\omega$},
$$
we obtain
$$
\sum_{k=1}^N Q_k(a-\www{a})t^{-\alpha m(k)}
+ \sum_{k=1}^N S_{k,\ell_1}(\mu_{\ell_1}f - \www{\mu}_{\ell_1}\www{f})
t^{-\alpha m(k) - \ell_1-1}
$$
$$
= O(t^{-\alpha m(1) - \ell_1 -2}) 
+ O(t^{-\alpha m(N+1)}) \quad \mbox{in $\omega$}      \eqno{(4.2)}
$$
as $t \to \infty$.

Thanks to (1.11), we can argue in the same way as Second Step of the proof
of Theorem 1 in Section 3, by (4.1) we can verify
$$
a = \www{a} \quad \mbox{and} \quad 
f = \frac{\www{\mu}_{\ell_1}}{\mu_{\ell_1}}\www{f} \quad 
\mbox{in $\OOO$.}
                                \eqno{(4.3)}
$$
Since $u$ and $\www{u}$ satisfies (1.1) and (1.8), setting
$y:= u - \www{u}$ in $\OOO \times (0,\infty)$, in view of (4.3)
we obtain
$$
\left\{ \begin{array}{rl}
& \ddda y = -Ay + \rho(t)\www{f}(x), \quad x\in \OOO,\, t>0, \\
& y\vert_{\ppp\OOO\times (0,\infty)} = 0, \\
& y(\cdot,0) = 0 \quad \mbox{in $\OOO$},
\end{array}\right.
$$
where 
$$
\rho(t) := \frac{\www{\mu}_{\ell_1}}{\mu_{\ell_1}}\mu(t) - \www{\mu}(t),
\quad t>0.
$$
Moreover by Lemma 4, we have
$$
y(x,t) = 0, \quad x\in \omega, \, t>T.
$$
Therefore, in view of (2.1), we apply (3.1) to have
$$
\sum_{k=1}^N  \frac{(-1)^{m(k)}}{\Gamma(-\alpha m(k))} 
\left( \int^T_0 (t-s)^{-\alpha m(k) -1} \rho(s) ds\right)
\sumn \frac{P_n\www{f}}{\la_n^{m(k)+1}}
$$
$$
= O( t^{-\alpha m(N+1)-1}) \quad \mbox{in $\omega$}             \eqno{(4.4)}
$$
as $t \to \infty$.

Assume that $\rho\not\equiv 0$ in $(0,T)$, that is,
$$
\www{\mu}_{\ell_1}\mu(t) \not\equiv \mu_{\ell_1}\www{\mu}(t) \quad 
\mbox{in $(0,T)$}.                                 \eqno{(4.5)}
$$
Then the Weierstrass theorem yields the existence of $\ell_2 \in 
\N \cup \{ 0\}$ such that 
$$
\rho_{\ell_2}: = \int^T_0 (-s)^{\ell_2} \rho(s) ds \ne 0, \quad
\int^T_0 (-s)^m \rho(s) ds = 0 \quad\mbox{for $0\le m \le \ell_2-1$
if $\ell_2 \in \N$}.
$$
Therefore, the same agument as in (3.3) implies
$$
\int^T_0 (t-s)^{-\alpha m(k)-1} \rho(s) ds
= \rho_{\ell_2}\left(
\begin{array}{c}
-\alpha m(k) - 1 \\
\ell_2 \\
\end{array}
\right) t^{-\alpha m(k)-1-\ell_2}\left( 1 + O\left(\frac{1}{t}
\right)\right)
$$
as $t \to \infty$.

Substituting this into (4.4), we obtain 
\begin{align*}
&\sum_{k=1}^N \rho_{\ell_2} \frac{(-1)^{m(k)}}{\Gamma(-\alpha m(k))}
\left(
\begin{array}{c}
-\alpha m(k) -1 \\
\ell_2 \\
\end{array}
\right) t^{-\alpha m(k) - 1 - \ell_2}
\left( \sumn \frac{P_n\www{f}}{\la_n^{m(k)+1}}\right) \\
= & O(t^{-\alpha m(1) - \ell_2 -2}) + O(t^{-\alpha m(N+1) - 1})
\end{align*}
in $\omega$ as $t \to \infty$, that is,
$$
\sum_{k=1}^N \rho_{\ell_2} S_{k,\ell_2}(\www{f})t^{-\alpha m(k) - 1-\ell_2}
= O(t^{-\alpha m(1) - \ell_2-2}) + O(t^{-\alpha m(N+1) - 1}) 
\quad \mbox{in $\omega$}
$$
as $t \to \infty$.

Arguing in the same way as in Second Step of the proof of Theorem 1,
we can see $\rho_{\ell_2} S_{k, \ell_2}(\www{f}) =  0$ in $\omega$ for 
each $k \in \N$.  The assumption of the theorem implies
$\rho_{\ell_2} \ne 0$.  Hence, $S_{k, \ell_2}(\www{f}) =  0$ in $\omega$ for 
each $k \in \N$, that is,
$\sumn \frac{P_n\www{f}}{\la_n^{m(k)+1}} = 0$ in $\omega$ for 
each $k\in \N$.
Therefore, in view of Lemma 6, applying the same argument at the end of 
Second Step of the proof of Theorem 1 in Section 3, we reach 
$\www{f} = 0$ in $\OOO$.  Then (4.3) implies $f=0$ in $\OOO$, which is 
a contradiction for $f\not\equiv 0$ in $\OOO$ or
$\www{f}\not\equiv 0$ in $\OOO$.
Hence assumption (4.5) is impossible, and so 
$\www{\mu}_{\ell_1}\mu(t) = \mu_{\ell_1}\www{\mu}(t)$ for 
$0 < t < T$.
Thus the proof of Theorem 2 (i) is complete.
$\blacksquare$
\\
{\bf Proof in the case $1<\alpha<2$.}
\\
We assume that $\mu_{\ell_1}:= \int^T_0 (-s)^{\ell_1} \mu(s) ds
\ne 0$.  Arguing the same way for (4.2), we obtain
\begin{align*}
& \sum_{k=1}^N Q_k(a-\www{a})t^{-\alpha m(k)}
+ \sum_{k-1}^N R_k(b-\www{b})t^{-\alpha m(k) + 1}\\
+ & \sum_{k=1}^N S_{k,\ell_1}(\mu_{\ell_1}f-\www{\mu}_{\ell_1}\www{f}) 
t^{-\alpha m(k) - \ell_1-1}
= O(t^{-\alpha m(1) - \ell_1-2}) + O(t^{-\alpha m(N+1)+1})
\end{align*}
in $\omega$ as $t \to \infty$.

Thus in the same way as in Third Step of the proof of Theoerm 1 in 
Section 3, we can obtain $a = \www{a}$, $b = \www{b}$ in $\OOO$ and 
$\mu_{\ell_1}f = \www{\mu}_{\ell_1}\www{f}$ in $\OOO$.

Now the same arguments as in the case $0<\alpha<1$ after (4.3) bring
$\www{\mu}_{\ell_1}\mu(t) = \mu_{\ell_1}\www{\mu}(t)$ for 
$0<t<T$.  Thus the proof of Theorem 2 (ii) is complete.
$\blacksquare$ 
  
\section{Concluding remarks}

{\bf 1.} 
We have proved the uniqueness in simultaneously determining 
factors $f(x)$ and $\mu(t)$ of a source term and initial values $a$ 
and/or $b$, except for some values of orders $\alpha$.
The proof is by equating the terms $t^{-\alpha m(k)}$, etc. 
with the same orders of the asymptotic expansions of the solutions.
Those terms are generated by initial data and a source term:  
$\int^T_0 (t-s)^{-\alpha m(k) - 1}\mu(s) ds$ with $k\in \N$ (see
(3.5) and (3.6)).
Without conditions (1.6), (1.7), (1.11) and (1.12) on $\alpha$, 
infinitely many terms are overlapped and our argument cannot conclude the 
uniqueness. So far, in general, we do not know the uniqueness without such 
conditions on $\alpha$.
However in a case of time-fractional ordinary differential equations,
although the separation of the terms is not perfect, we can prove the 
uniqueness corresponding to Theorem 2, 
thanks to that an unknown initial value is scalar.
More precisely,
\\
{\bf Proposition 1.}
\\
{\it
Let $\alpha \in (0,1) \cup (1,2)$ be arbitrarily given.
For $\la \in \R$, we consider
$$
\left\{ \begin{array}{rl}
& d_t^{\alpha} u(t) = -\la u(t) + \mu(t), \\
& u(0) = a, \quad \mbox{if $0<\alpha<1$}, \\
& u(0) = a, \quad \frac{du}{dt}(0) = b, \quad \mbox{if $1<\alpha<2$}.
\end{array}\right.
                               \eqno{(5.1)}
$$
We assume that $\mu(t)$ satisfies (1.4), and $T < T_1 < T_2$.  Then 
\\
(i) Case $0<\alpha<1$: if $u(t) = 0$ for $T_1 < t < T_2$, then 
$a=0$ and $\mu(t) = 0$ for $0<t<T$.
\\
(ii) Case $1<\alpha<2$: if $u(t) = 0$ for $T_1 < t < T_2$, then 
$a=b=0$ and $\mu(t) = 0$ for $0<t<T$.
}
\\
{\bf Proof of Proposition 1.}
\\
It is known (e.g., Kilbas, Srivastava and Trujillo \cite{KST}, p.141) that 
$$
u(t) =
\left\{ \begin{array}{rl}
& aE_{\alpha,1}(-\la t^{\alpha}) + \int^t_0 (t-s)^{\alpha-1}
E_{\alpha,\alpha}(-\la (t-s)^{\alpha}) \mu(s) ds, \quad 0<\alpha<1, \\
& aE_{\alpha,1}(-\la t^{\alpha}) + bt E_{\alpha,2}(-\la t^{\alpha}) \\
+ &\int^t_0 (t-s)^{\alpha-1}
E_{\alpha,\alpha}(-\la (t-s)^{\alpha}) \mu(s) ds, \quad 1<\alpha<2,
\quad t>0.
\end{array}\right.
                                            \eqno{(5.2)}
$$
In view of (3.1) with $N=1$, we have
$$
\left\{ \begin{array}{rl}
& E_{\alpha,1}(-\la t^{\alpha}) 
= \frac{(-1)^{m(1)+1}}{\Gamma(1-\alpha m(1))}\frac{1}{\la^{m(1)}}
\frac{1}{t^{\alpha m(1)}} + O(t^{-\alpha m(2)}), \\
& tE_{\alpha,2}(-\la t^{\alpha}) 
= \frac{(-1)^{m(1)+1}}{\Gamma(2-\alpha m(1))}\frac{1}{\la^{m(1)}}
\frac{1}{t^{\alpha m(1)-1}} + O(t^{-\alpha m(2)+1}), \\
& t^{\alpha-1}E_{\alpha,\alpha}(-\la t^{\alpha}) 
= \frac{(-1)^{m(1)}}{\Gamma(-\alpha m(1))}\frac{1}{\la^{m(1)+1}}
\frac{1}{t^{\alpha m(1)+1}} + O(t^{-\alpha m(2)-1})
\end{array}\right.
                                         \eqno{(5.3)}
$$
as $t \to \infty$.
    
Since $u(t) = 0$ for $T_1<t<T_2$, similarly to Lemma 4, we obtain
$$
u(t) = 0, \quad t> T.
$$
\\
{\bf Case: $0<\alpha<1$.}
\\
By (5.2) and $\mu(t) = 0$ for $t>T$, we see
$$
aE_{\alpha,1}(-\la t^{\alpha}) + \int^T_0 (t-s)^{\alpha-1}
E_{\alpha,\alpha}(-\la (t-s)^{\alpha}) \mu(s) ds = 0, \quad t>2T.
                                                      \eqno{(5.4)}
$$
Therefore, by (5.3), we obtain
\begin{align*}
&\frac{a(-1)^{m(1)+1}}{\Gamma(1-\alpha m(1))}\frac{1}{\la^{m(1)}}
\frac{1}{t^{\alpha m(1)}} + O(t^{-\alpha m(2)})\\
+ &\frac{(-1)^{m(1)}}{\Gamma(-\alpha m(1))}\frac{1}{\la^{m(1)+1}}
\int^T_0 (t-s)^{-\alpha m(1)-1} \mu(s) ds
+ \int^T_0 O((t-s)^{-\alpha m(2)-1}) \mu(s) ds,
\end{align*}
and so the multiplication by $t^{\alpha m(1)}$ yields
$$
\frac{a(-1)^{m(1)+1}}{\Gamma(1-\alpha m(1))}\frac{1}{\la^{m(1)}}
+ \frac{(-1)^{m(1)}}{\Gamma(-\alpha m(1))}\frac{1}{\la^{m(1)+1}}
\int^T_0 \frac{t^{\alpha m(1)}}{(t-s)^{\alpha m(1)+1}} \mu(s) ds
$$
$$
+ O(t^{-\alpha (m(2)-m(1))})
+ \int^T_0 O(t^{\alpha m(1)}(t-s)^{-\alpha m(2)-1}) ds
= 0 \quad \mbox{as $t \to \infty$}.
                                     \eqno{(5.5)}
$$
Here, for $0 < s < T < 2T < t$, we have
$\frac{1}{1-\frac{s}{t}}\le 2$ and $\frac{1}{t-s} \le \frac{1}{T}$, so that 
$$
\frac{t^{\alpha m(1)}}{(t-s)^{\alpha m(1)+1}} 
= \left( \frac{1}{1 - \frac{s}{t}}\right)^{\alpha m(1)} \frac{1}{t-s}
\le 2^{\alpha m(1)}\frac{1}{T} \quad \mbox{for all $s< T < 2T < t$}.
$$
Since $\lim_{t\to \infty} \frac{t^{\alpha m(1)}}{(t-s)^{\alpha m(1)+1}} = 0$,
the Lebesgue convergence theorem implies 
$$
\lim_{t\to \infty} \int^T_0 \frac{t^{\alpha m(1)}}{(t-s)^{\alpha m(1)+1}} 
\mu(s) ds = 0.                            \eqno{(5.6)}
$$
Therefore, letting $t \to \infty$ in (5.5), we see
$$
\frac{a(-1)^{m(1)+1}}{\Gamma(1-\alpha m(1))}\frac{1}{\la^{m(1)}} = 0,
$$
which implies $a=0$.  Hence (5.4) yields
$$
\int^T_0 (t-s)^{\alpha-1} E_{\alpha,\alpha}(-\la (t-s)^{\alpha}) \mu(s) ds
= 0, \quad t>2T.                        \eqno{(5.7)}
$$

Therefore, the following Lemma 7 yields that $\mu = 0$ in $(0,T)$, and so 
completes the proof of Proposition 1 for the case $0<\alpha<1$.
Thus it suffices to prove
\\
{\bf Lemma 7.}
\\
{\it
Let $\mu \in L^{\infty}(0,T)$.  Then (5.7) implies 
$\mu = 0$ in $(0,T)$.
}
\\
{\bf Proof of Lemma 7.}
Substituting the third asymptotic expansion in (3.1) into (5.7), we have
\begin{align*}
& 0 = \int^T_0 (t-s)^{\alpha-1}E_{\alpha,\alpha}(-\la(t-s)^{\alpha}) \mu(s) ds
                     \\
= &\sum_{k=1}^N \frac{(-1)^{m(k)}}{\Gamma(-\alpha m(k))\la^{m(k)+1}}
\left( \int^T_0 (t-s)^{-\alpha m(k) -1} \mu(s) ds \right)
+ O(t^{-\alpha m(N+1) - 1})
\end{align*}
as $t \to \infty$ for all $N\in \N$.

Contrarily assume that $\mu \not\equiv 0$ in $(0,T)$.  Then Lemma 2 implies
the existence of $\ell_0 \in \N \cup \{ 0\}$ such that 
$$
\mu_{\ell_0}:= \int^T_0 (-s)^{\ell_0} \mu(s) ds, \quad
\int^T_0 (-s)^m \mu(s) ds \quad \mbox{for $0 \le m \le \ell_0-1$ if
$\ell_0\in \N$.}
$$
Applying (3.4), we obtain
$$
\sum_{k=1}^N \frac{(-1)^{m(k)}}{\Gamma(-\alpha m(k)) \la^{m(k)+1}}
t^{-\alpha m(k) - \ell_0-1} \mu_{\ell_0}
\left(
\begin{array}{c}
-\alpha m(k) - 1 \\
\ell_0 \\
\end{array}
\right) 
\left( 1 + O\left( \frac{1}{t}\right)\right)
+ O(t^{-\alpha m(N+1) -1}) = 0
$$
as $t \to \infty$.  Lemma 5 yields
$$
\mu_{\ell_0}\frac{(-1)^{m(1)}}{\Gamma(-\alpha m(1)) \la^{m(1)+1}}
\left(
\begin{array}{c}
-\alpha m(1) - 1 \\
\ell_0 \\
\end{array}
\right) = 0.
$$
Since 
$$
\frac{(-1)^{m(1)}}{\Gamma(-\alpha m(1)) \la^{m(1)+1}}
\left(
\begin{array}{c}
-\alpha m(1) - 1 \\
\ell_0 \\
\end{array}
\right) \ne 0
$$
by $\alpha m(1)\not\in \N \cup \{ 0\}$, we obtain
$\mu_{\ell_0} = 0$.  This is a contradiction by the definition of 
$\mu_{\ell_0}$.  Therefore $\int^T_0 (-s)^m \mu(s) ds = 0$ for all
$m\in \N \cup \{0\}$.  Lemma 2 yields that $\mu = 0$ in $(0,T)$.
Thus the proof of Lemma 7 is complete.
$\blacksquare$
\\
{\bf Case: $1<\alpha<2$.}
\\
By (5.2), we have
$$
aE_{\alpha,1}(-\la t^{\alpha}) + bt E_{\alpha,2}(-\la t^{\alpha}) \\
+ \int^T_0 (t-s)^{\alpha-1}
E_{\alpha,\alpha}(-\la (t-s)^{\alpha}) \mu(s) ds = 0, \quad t>2T.
                                            \eqno{(5.8)}
$$
Substituting (5.3) into (5.8), we obtain
$$
\frac{a(-1)^{m(1)+1}}{\Gamma(1-\alpha m(1))}\frac{1}{\la^{m(1)}}
\frac{1}{t^{\alpha m(1)}} 
+ \frac{b(-1)^{m(1)+1}}{\Gamma(2-\alpha m(1))}\frac{1}{\la^{m(1)}}
\frac{1}{t^{\alpha m(1)-1}} 
$$
$$
+ \frac{(-1)^{m(1)}}{\Gamma(-\alpha m(1))}\frac{1}{\la^{m(1)+1}}
\int^T_0 (t-s)^{-\alpha m(1)-1} \mu(s) ds
+ O(t^{-\alpha m(2)+1}) = 0 \quad \mbox{as $t \to \infty$}.
                                                     \eqno{(5.9)}
$$
Similarly to (5.6), we have
$$
\lim_{t\to \infty} \int^T_0 \frac{t^{\alpha m(1)-1}}{(t-s)^{\alpha m(1)+1}}
\mu(s) ds = 0.
$$
Consequently, multiplying by $t^{\alpha m(1)-1}$ and letting $t \to \infty$, 
we see
$$
b\frac{(-1)^{m(1)+1}}{\Gamma(2-\alpha m(1))}\frac{1}{\la^{m(1)}} = 0,
$$
that is, $b=0$.  Hence, multiplying by $t^{\alpha m(1)}$ and using (5.6) 
in (5.9), we obtain $a=0$.
Therefore (5.8) implies (5.7).  Lemma 7 yields $\mu(t) = 0$ for $0<t<T$.
Thus the proof of Proposition 1 is complete.
$\blacksquare$

{\bf 2.}
We take interior data $u\vert_{\omega\times (T_1,T_2)}$ as 
observation data.
We can discuss to obtain the uniqueness with other typs of data such as
boundary data $\ppp_{\nu}u\vert_{\gamma \times (T_1,T_2)}$, where 
$\gamma\subset \ppp\OOO$ is an arbitrarily chosen subboundary.

{\bf 3.}
We should extend the uniqueness to more general elliptic 
operator, especially, non-symmetric $A$.   
As for the determination of initial values for $\alpha \in (1,2)$, 
see Loreti, Sforza and Yamamoto \cite{LSY}.
The extension of the uniqueness in the current article to non-symmetric 
$A$ should be a future work.

{\bf 4.}
Our argument is applicable to other type of inverse problem.
For example, for $\alpha \in (0,1) \cup (1,2)$, we consider 
$$
\left\{ \begin{array}{rl}
& \ddda y = -Ay + \mu(t)f(x), \quad x\in \OOO, \, t>0,\\
& y\vert_{\ppp\OOO \times (0,\infty)} = 0, \\
& y(\cdot, 0) = 0 \quad \mbox{in $\OOO$ $\quad$ if $0<\alpha < 1$},\\
& y(\cdot, 0) = \ppp_ty(\cdot,0) = 0 \quad \mbox{in $\OOO$ $\quad$ if 
$1<\alpha < 2$}.
\end{array}\right.
$$
We assume that $\mu$ satisfies (1.4).  Then we dicuss
\\
{\bf Inverse source problem.}
{\it We choose $x_0 \in \OOO$ and $T_1, T_2 > 0$ such that $T<T_1<T_2$.
Then determine $\mu(t)$ for $0<t<T$ by data $y(x_0,t)$ for $T_1<t<T_2$.
}

For simplicity, let $f \in C^{\infty}_0(\OOO)$. 
We can prove
\\
{\bf Proposition 2.}
\\
{\it 
We assume that $f(x_0) \ne 0$.  Then 
$y(x_0,t) = 0$ for $T_1 < t < T_2$ implies 
$\mu(t) =  0$ for $0<t<T$.
}

The uniqueness for $\alpha \in (0,1) \cup (1,2)$ is different from the 
cases $\alpha =1$ and $=2$.  In particular, in the case where $\alpha=1$ and
the spatial dimension is one, we have no uniqueness (\cite{CLY}).
\\
{\bf Proof of Proposition 2.}
The proof is based on (4.4).
Assume that $\mu \not\equiv 0$ in $(0,T)$.  Then Lemma 2 
implies that there exists $\ell_0 \in \N \cup \{ 0\}$ satisfying 
(1.5), and in particular, we have $\mu_{\ell_0}:= \int^T_0 (-s)^{\ell_0}
\mu(s) ds \ne 0$.
By $f\in C^{\infty}_0(\OOO)$, we can verify that 
$\sumn P_nf$ converges in $C(\ooo{\OOO})$.  Therefore, similarly to (4.4), 
we can obtain
\begin{align*}
& y(x_0,t) = \sum_{k=1}^N  \frac{(-1)^{m(k)}}{\Gamma(-\alpha m(k))} 
\left( \int^T_0 (t-s)^{-\alpha m(k) -1} \mu(s) ds\right)
\sumn \frac{(P_nf)(x_0)}{\la_n^{m(k)+1}}                 \\
= & O( t^{-\alpha m(N+1)-1})
\end{align*}
for all $t > 2T$ and $N\in \N$.

Applying (3.4), we reach 
\begin{align*}
& \sum_{k=1}^N  \frac{(-1)^{m(k)}}{\Gamma(-\alpha m(k))} 
\mu_{\ell_0}
\left(
\begin{array}{c}
-\alpha m(k) - 1 \\
\ell_0 \\
\end{array}
\right) t^{-\alpha m(k) - \ell_0 - 1}
\left( 1 + O\left( \frac{1}{t}\right)\right)
\sumn \frac{(P_nf)(x_0)}{\la_n^{m(k)+1}}\\
=& O( t^{-\alpha m(N+1)-1}) \quad \mbox{for all 
$t > 2T$ and all $N\in \N$}.
\end{align*}
In view of Lemma 5, we can prove
$$
\sumn \frac{(P_nf)(x_0)}{\la_n^{m(k)+1}} = 0 \quad 
\mbox{for all $k\in \N$}.
$$
Hence Lemma 6 yields $(P_nf)(x_0) = 0$ for all $n\in \N$.
In terms of $f \in C^{\infty}_0(\OOO)$, we see that 
$f(x_0) = \sumn (P_nf)(x_0) = 0$, which contradicts the 
assumption $f(x_0) \ne 0$.  Therefore, $\mu \equiv 0$ in $(0,T)$.
Thus the proof of Proposition 2 is complete.
$\blacksquare$
\\

{\bf Acknowledgements.} 

The third author is supported by Grant-in-Aid for Scientific Research (A)
20H00117 and Grant-in-Aid for Challenging Research (Pioneering) 21K18142, 
Japan Society for the Promotion of Science (JSPS).


\end{document}